\documentclass[reqno,11pt,a4paper]{article}
\usepackage{amssymb,amsthm,amsmath,dsfont,graphicx,fullpage,ccaption}
\usepackage[usenames,dvipsnames]{color}

\newfixedcaption{\figcaption}{figure}
\newfixedcaption{\tablecaption}{table}

\ifx\pdfoutput\undefined
\usepackage[dvips,bookmarks]{hyperref}
\else
\usepackage{hyperref}
\fi
\hypersetup{%
pdftitle={M/G/oo polling systems with random visit times},
pdfauthor={M. Vlasiou, U. Yechiali},
pdfkeywords={},
bookmarks=true,%
bookmarksnumbered=true,
pdfstartview={FitH},
colorlinks=true,
citecolor=Plum,
linkcolor=RedOrange,
urlcolor=RawSienna,
bookmarksopen=true,%
bookmarksopenlevel=1,
unicode=true,%
breaklinks=true,%
}%

\theoremstyle{plain}
\newtheorem{theorem}{Theorem}

\theoremstyle{remark}
\newtheorem{remark}{Remark}

\numberwithin{equation}{section}

\newcommand{\e}{\ensuremath{\mathbb{E}}}
\newcommand{\p}{\ensuremath{\mathbb{P}}}
\renewcommand{\d}{\ensuremath{\,\mathrm{d}}}

\newcommand{\ltd}[1][]{\ensuremath{\widetilde{D}_{#1}}}

\graphicspath{{graphs/}}

\begin{document}
\title{M/G/$\infty$ polling systems with random visit times}
\author{M.\ Vlasiou$^{*}$, U.\ Yechiali$^{**}$}
\date{June 7, 2007}
\maketitle

\begin{center}
$^{*}$ Georgia Institute of Technology,
\\H.\ Milton Stewart School of Industrial \& Systems Engineering,
\\765 Ferst Drive, Atlanta GA 30332-0205, USA.
\end{center}

\begin{center}
$^{**}$ Department of Statistics and Operations Research,
\\School of Mathematical Sciences,
\\Raymond and Beverly Sackler Faculty of Exact Sciences,
\\Tel Aviv University, Tel Aviv 69978, Israel.
\\\vspace{0.3cm}\href{mailto:vlasiou@gatech.edu}{vlasiou@gatech.edu}, \href{mailto:uriy@post.tau.ac.il}{uriy@post.tau.ac.il}
\end{center}


\begin{abstract}
We consider a polling system where a group of an infinite number of servers visits sequentially a set of queues. When visited, each queue is attended for a random time. Arrivals at each queue follow a Poisson process, and service time of each individual customer is drawn from a general probability distribution function. Thus, each of the queues comprising the system is, in isolation, an M/G/$\infty$-type queue.  A job that is not completed during a visit will have a new service time requirement sampled from the service-time distribution of the corresponding queue. To the best of our knowledge, this paper is the first in which an M/G/$\infty$-type polling system is analysed. For this polling model, we derive the probability generating function and expected value of the queue lengths, and the Laplace-Stieltjes transform and expected value of the sojourn time of a customer. Moreover, we identify the policy that maximises the throughput of the system per cycle and conclude that under the Hamiltonian-tour approach, the optimal visiting order is \emph{independent} of the number of customers present at the various queues at the start of the cycle.
\end{abstract}

\section{Introduction}\label{s:intro}
A typical polling system consists of a number of queues, attended by a single server in a cyclic fashion. There is a huge body of literature on polling systems that has developed since the late 1950s, when the papers of Mack et al.~\cite{mack57a,mack57} concerning a patrolling repairman model for the British cotton industry were published. Rather than giving a partial overview of the literature, we refer the interested reader to the following books, surveys, and papers on polling systems: Takagi~\cite{takagi-APS,takagi90,takagi97}, Boxma and Groenendijk~\cite{boxma87}, Levy and Sidi~\cite{levy90}, Yechiali~\cite{yechiali93}, Borst~\cite{borst}, Eliazar and Yechiali~\cite{eliazar98}, Nakdimon and Yechiali~\cite{nakdimon03}.

Polling systems have been used as a central model for the analysis of a wide variety of applications in the areas of repair problems~\cite{mack57a,mack57}, telecommunication systems~\cite{cooper69}, road traffic control~\cite{stidham69}, computer networks~\cite{wang85}, multiple access protocols~\cite{bernabeu95},  multiplexing schemes in ISDN~\cite{twu96},  satellite systems~\cite{altman02}, flexible manufacturing systems~\cite{vuuren06}, and the like.

In many of these applications, as well as in most polling models, it is customary to control the amount of service given to each queue during the server's visit. Common service policies are the \textit{exhaustive}, \textit{gated}, \textit{globally gated} and \textit{limited} regimes. Under the exhaustive regime, at each visit the server attends the queue until it becomes completely empty, and only then is the server allowed to move on. Under the gated regime, the only customers served during a visit are the ones who are present when the server enters (polls) the queue, while customers arriving when the queue is attended will be served during the next visit. The globally gated regime, introduced by Boxma, Levy and Yechiali~\cite{boxma92}, is a modification of the gated one: the only customers served during a visit are those who are present at the beginning of a cycle. Finally, under the $k$-Limited service discipline only a limited number of jobs (at most $k$) are served at each server's visit to each queue. These service policies imply that the duration of the visit time in a polled queue is a function of the number of customers present there at a given moment (such as the beginning of the cycle or the moment the server enters the queue).

In this paper, we analyse a polling system that differs in two ways from the classical polling model. Rather than considering a \emph{single} server providing service to customers at the various queues, we assume that an \emph{infinite number} of servers is moving as a single group between the queues. Moreover, the service policy we study is \emph{independent} of the queue length. We assume that the group of servers visits each queue for a (possibly random) amount of time that is independent of everything else and which has a distribution that may vary per queue. We further assume that the arrival process of customers to each queue is Poisson and that the service time distribution for customers in each queue is general. To the best of our knowledge, this paper is the first in which an M/G/$\infty$-type polling system is analysed.

The specific application that raised our attention and led us to this model is in the field of road traffic control. Polling models for road traffic are typically along the lines of the classical polling system; namely, they involve a single server rotating around a number of queues. Other assumptions that are typically being made for such models include deterministic service times (i.e.\ the amount of time that a car needs to pass a traffic light after possibly standing in the queue) and deterministic visit times (i.e.\ the time the traffic light remains green); see, for example, van der Heijden~\cite{heijden01}. Although these models provide fairly good approximations of reality, such assumptions fail to capture the variation both in service times and in visit times. Cars do not need the same amount of time to cross a segment of the road; the ones standing ahead in the queue will inevitably need less time and those that arrive while the queue is empty and the traffic light is still green will not even require the additional time incurred by acceleration. Moreover, recent developments in the technology of traffic lights has led to the design of traffic lights that do not turn green unless a queue is formed, and turn red either when the queue is empty or after a maximum amount of time, which may also vary within a day. As a result, in this paper we provide a framework for studying road traffic control under less restrictive assumptions. We propose an infinite-server polling system, which models the behaviour of traffic: while the traffic light is green all cars present in the queue or approaching the traffic light proceed (receive service) and the time they need to complete service is assumed to be a random variable following a general distribution. Furthermore, we assume that the time the traffic light is green (visit time) is random, although our results are directly applicable in case of deterministic visit times or, more generally, in case the visit times follow a discrete distribution taking positive values.

A common approximation to road traffic is to consider the traffic as \emph{fluid} passing through the road. This approximation is fairly accurate when the traffic is relatively high. Mathematically, high traffic can be modelled by assuming that the arrival rate of customers at each of the queues tends to infinity. The study of such a model provides insights at the queue length (and thus the congestion of a junction) under heavy load. In this paper though we do not study the evolution of the system under heavy load. We assume that the arrival rate at each queue is fixed. This assumption is usually made for the standard polling systems and provides a reasonable approximation to normal traffic conditions.

The rest of the present paper is organized as follows. Section~\ref{s:model} introduces the model, gives further notation, and describes formally the evolution of the system. In Section~\ref{s:queue} we compute recursively the first moment and the probability generating function  of the queue length distributions at a polling instant. Later on, in Section~\ref{s:sojourn} we derive the mean and the Laplace-Stieltjes transform of the sojourn time of a customer arriving at queue $i$, and we show how these expressions simplify in the special case where both the service time and the visit time at queue $i$ are exponentially distributed. Based on the results derived up to that point, in Section~\ref{s:Numerical results} we give some numerical results. Specifically, we examine numerically the effect of the first two moments of the visit and service times on the sojourn time of an arbitrary customer. These numerical results indicate that there is an optimal value for the mean visit time to the various queues that minimises the mean sojourn time of an arbitrary customer. In Section~\ref{s:optimisation} we investigate how we can optimise the visit order of the servers at the various queues so that the expected throughput of the system is maximised. It emerges that even when considering \emph{semi-dynamic} control policies, in which the group of servers plans a new route for each cycle, the optimal visiting order that maximises the expected throughput per cycle is \emph{fixed} for all cycles. In other words, because of the infinite number of servers, information regarding the queue lengths of all queues at the beginning of a cycle has no effect on the choice of the optimal strategy.

\section{Model description and notation}\label{s:model}
We consider a polling system with $N \geqslant 2$ infinite-buffer queues attended by a group of ample number of servers that visits the queues in a fixed cyclic fashion. We index the queues by $i=1,2,\ldots,N$ in the order of the servers' movement. We shall refer to the polling instant of queue $i$ as the moment when the servers enter that queue. When visiting queue $i$, the group of servers continues working at this queue for $V_i$ units of time, and acts there as an M/G/$\infty$ queue. We assume that the visit times are independent, identically distributed (i.i.d.)\ random variables.

Customers arrive at all queues according to independent homogeneous Poisson processes with rate $\lambda_i$ for queue $i$. After completing their service time, customers leave the system. The service time of each individual customer at queue $i$ is denoted by $B_i$. It is assumed that all service times in one queue are i.i.d.\ random variables, which are mutually independent of all service times at any other queue. At the end of a visit to queue $i$, the group of servers moves to queue $i+1$, incurring a switch-over time $D_i$ and a realisation of $V_{i+1}$ is drawn. We assume that $\{D_i\}$ is a sequence of independent random variables. The total switch-over time during a full cycle is $D=\sum_{i=1}^{N} D_i$, and the length of a full cycle is denoted by the random variable $C$. We assume that all random variables so far are mutually independent.

During the visit time of the group of servers to queue $i$, a customer present at queue $i$ at the polling instant of that queue will successfully complete his service with probability {$p_i(V_i)=\p[B_i\leqslant V_i \mid V_i]$}. We assume that if the service of a customer of queue $i$ is not completed during a single visit, then at the next visit a new service time will be drawn from the service time distribution of $B_i$ for that particular customer.

For a generic random variable $Y_i$, we denote its first two moments by $\e[Y_i]$ and $\e[Y_i^2]$, respectively. Thus, for example, $\e[V_i]$ is the mean visit time of the servers at queue $i$. By convention, $\sum_{i\neq j} Y_i=\sum_{\stackrel{i=1}{i\neq j}}^N Y_i$, and similarly for the product operator.  All further notation will be introduced when it is first used.

\subsection*{Law of motion}
Let $X_i^j$, $i,j=1,2,\ldots,N$, denote the number of customers in queue $j$ at the moment when queue $i$ is polled {and let $A_j(t)$ denote the number of Poisson arrivals to queue $j$ during a time interval of length $t$}. The law of motion describing the evolution of the system when the server moves from queue $i$ to queue $i+1$ connects  $X_{i+1}^j$ to $X_i^j$ and is given by
\begin{equation}\label{eq:lom}
  X_{i+1}^j=\begin{cases}
                X_i^j+{A_{j,1}(V_i)+A_{j,2}(D_i)},     &j\neq i,\\
                \mathrm{Binom}(X_i^i, 1 - {p_i(V_i)})+\mathrm{Poisson}\big(\Lambda_i(V_i)\big)+A_i(D_i), &j=i,
  \end{cases}
\end{equation}
where {for all $k$, $A_{j,k}(t)$ is an i.i.d.\ copy of $A_j(t)$}, Binom$(n,p)$ is a binomial random variable with parameters $n$ and $p$, and Poisson$\big(\Lambda_i(t)\big)$ is a Poisson random variable with rate
$$
\Lambda_i(t)=\lambda_i\int_0^t \p[B_i>y] \d y.
$$

{Note that from \eqref{eq:lom} we see that for all $j$, the random variables $X_i^j$ are independent of $V_i$ and $D_i$, which is evident, considering that the number of customers in a queue at the beginning of a visit does not depend on the length of the upcoming visit time or switch-over time.}

The relation for $j\neq i$ is straightforward. The number of customers at queue $j$ at polling instant of queue $i+1$ equals the number of customers that were there at polling instant of queue $i$ plus all customers that arrived during the visit time of queue $i$ and the switch-over time from queue $i$ to queue $i+1$.

For $j=i$, the relation is more involved. When the servers start polling queue $i$ they encounter $X_i^i$ customers. After $V_i$ time units, only a binomial number of customers out of the initial $X_i^i$ is still present. The probability that a single customer does not complete his service after $V_i$ time units is $1-p_i(V_i)=\p[B_i>V_i \mid V_i]$. In addition, there is a stream of new arrivals to queue $i$. The number of customers present at time $t$ in an M/G/$\infty$ queue (starting with zero customers at time $t=0$) is Poisson distributed with rate $\Lambda_i(t)$, as it is given above; see Tak\'{a}cs~\cite{takacs-ITQ}. The last term at the right-hand side of \eqref{eq:lom} incorporates the customers that arrived at the queue during the switch-over time from queue $i$ to queue $i+1$.

We shall employ this relation to derive the mean queue length and the probability generating function of all queues at a polling instant.

\section{Queue lengths at polling instants}\label{s:queue}
One of the main tools used in the analysis of polling systems is the derivation of a set of multi-dimensional probability generating functions of the number of jobs present in the various queues at a polling instant of queue $i$. The common method is to derive the probability generating function of a given queue at some polling instant in terms of the probability generating function of the same queue at the previous polling instant. Then, from the set of $N$ (implicit) dependent equations of the unknown probability generating functions one can obtain expressions which allow for numerical calculation of the mean queue length at each queue. These equations simplify significantly for several cases of the distribution of the visit times. In this section, we use the law of motion (buffer occupancy), which is given by Equation~\eqref{eq:lom} and apply this technique to compute recursively the first moment and the probability generating function  of the queue length distributions at a polling instant.

\subsection{Mean queue length}
From \eqref{eq:lom} we have the following relation for the mean queue length of queue $j$ at two consecutive polling instants.
\begin{equation}\label{eq:means rec}
  \e[X_{i+1}^j]=\begin{cases}
    \e[X_i^j]+\lambda_j \e[V_i]+\lambda_j \e[D_i],  &j\neq i,\\
    (1-p_i)\e[X_i^i] +\e[\Lambda_i(V_i)] +\lambda_i \e[D_i], &j=i,
  \end{cases}
\end{equation}
where $p_i=\p[B_i \leqslant V_i]=\e[p_i(V_i)]$. Summing \eqref{eq:means rec} over $i$ we obtain
\begin{equation}\label{eq:jj means}
p_j\e[X^j_j] =\lambda_j \sum_{i\neq j} \e[V_i]+\e[\Lambda_j(V_j)]+\lambda_j \e[D].
\end{equation}
Indeed, in steady state, the mean number of jobs in queue $j$ at a polling instant equals the fraction of jobs $(1-p_j)\e[X^j_j] $ left behind at the end of the previous visit, plus the mean number of arrivals during the cycle time out of queue $j$, which is $\lambda_j \Bigl(\sum_{i\neq j} \e[V_i]+\e[D]\Bigr)$, plus the mean number of customers in a M/G/$\infty$ queue at time $V_j$. The mean queue length of queue $j$ at polling instant of queue $i$ is easily derived from \eqref{eq:means rec}, yielding
\begin{equation}\label{eq:means}
\e[X_i^j]=\e[X_j^j](1-p_j)+\e[\Lambda_j(V_j)]+\lambda_j \sum_{k=j+1}^{i-1} \e[V_k]+\lambda_j \sum_{k=j}^{i-1} \e[D_k].
\end{equation}

For example, suppose that $B_j$ is exponentially distributed with parameter $\mu_j$. Then,
$$
\Lambda_j(V_j)=\lambda_j\int_0^{V_j} e^{-\mu_j y} \d y=\frac{\lambda_j}{\mu_j}(1-e^{-\mu_j V_j}).
$$
Thus, $\e[\Lambda_j(V_j)]=\lambda_j(1-\e[e^{-\mu_j V_j}])/{\mu_j}$. So, in particular, if $V_j$ is also exponentially distributed with parameter $\gamma_j$, then we have that $\e[\Lambda_j(V_j)]=\lambda_j/(\gamma_j+\mu_j)$, and the mean queue length of each queue can now easily be computed recursively from \eqref{eq:means}.

\subsection{Recursive relation for the generating function}
Define the generating function of the queue length of all queues at polling instants of queue $i$ as $G_i (\mathbf{z})=\e[\prod_{j=1}^N z_j^{X_i^j}]$. Then, from \eqref{eq:lom} we have that
\begin{equation}\label{eq:gen fun rec}
G_{i+1}(\mathbf{z})=\e[\prod_{j\neq i} z_j^{X_i^j+A_{j,1}(V_i)+A_{j,2}(D_i)} z_i^{\mathrm{Binom}(X_i^i, 1-{p_i(V_i)})+\mathrm{Poisson}\big(\Lambda_i(V_i)\big)+A_i(D_i)}]
\end{equation}
By conditioning on the vector $(X_i^1,\ldots,X_i^N)$, on $V_i$, and on $D_i$, Equation \eqref{eq:gen fun rec} becomes
\begin{multline}\label{eq:gen fun broken}
G_{i+1}(\mathbf{z})=\e[\prod_{j\neq i} z_j^{X_i^j}\ \e[\prod_{j\neq i} z_j^{A_{j,1}(V_i)} \mid V_i]\ \e[\prod_{j=1}^N  z_j^{A_{j,2}(D_i)} \mid D_i] \times \\ \times {\e[z_i^{\mathrm{Binom}(X_i^i, 1-p_i(V_i))} \mid X^i_i,V_i]}\ \e[z_i^{\mathrm{Poisson}\big(\Lambda_i(V_i)\big)} \mid V_i]].
\end{multline}
Since the number of arrivals at any queue during a fixed amount of time is independent of the number of arrivals at any other queue during the same given period, we have that
\begin{align*}
\e[\prod_{j=1}^N  z_j^{A_j(D_i)} \mid D_i=x]&=\e[\prod_{j=1}^N  z_j^{A_j(x)}]=\prod_{j=1}^N\e[z_j^{A_j(x)}]\\
&=\prod_{j=1}^N \sum_{n=0}^\infty z_j^n \frac{(\lambda_j x)^n}{n!}\, e^{-\lambda_j x}=\prod_{j=1}^N e^{-\lambda_j x (1-z_j)}.
\end{align*}
Therefore, we have
\begin{equation}\label{eq:1}
\e[\prod_{j=1}^N  z_j^{A_j(D_i)} \mid D_i]=e^{-D_i \sum_{j=1}^N \lambda_j (1-z_j)}.
\end{equation}
Likewise, we obtain that
\begin{equation}
\e[\prod_{j\neq i} z_j^{A_j(V_i)} \mid V_i]=e^{-V_i \sum_{j\neq i} \lambda_j (1-z_j)}.
\end{equation}
Moreover,
\begin{align*}
\e[z_i^{\mathrm{Binom}(X_i^i, 1-p_i(V_i))} \mid X_i^i=k, V_i=x]&=\e[z_i^{\mathrm{Binom}(k, 1-p_i(x))}]\\
    &=\sum_{\ell=0}^{k} z_i^\ell \binom{k}{\ell}(1-p_i(x))^\ell p_i(x)^{k-\ell}=\big(p_i(x)+z_i[1-p_i(x)]\big)^{k},
\end{align*}
or in other words,
\begin{equation}
\e[z_i^{\mathrm{Binom}(X_i^i, 1-p_i(V_i))} \mid X_i^i, V_i]=\big(p_i(V_i)+z_i[1-p_i(V_i)]\big)^{X_i^i}.
\end{equation}
For the last term of the right-hand side of \eqref{eq:gen fun broken} we have that
\begin{align*}
\e[ z_i^{\mathrm{Poisson}\big(\Lambda_i(V_i)\big)}\mid V_i=x]&=\e[ z_i^{\mathrm{Poisson}\big(\Lambda_i(x)\big)}]=\sum_{n=0}^\infty z_i^n \frac{\big(\Lambda_i(x)\big)^n}{n!}\, e^{-\Lambda_i(x)}=e^{-\Lambda_i(x)(1-z_i)},
\end{align*}
which yields that
\begin{equation}\label{eq:3}
\e[ z_i^{\mathrm{Poisson}\big(\Lambda_i(V_i)\big)}\mid V_i]=e^{-\Lambda_i(V_i)(1-z_i)}.
\end{equation}
Substituting \eqref{eq:1} -- \eqref{eq:3} into \eqref{eq:gen fun broken}, we obtain
\begin{equation}\label{eq:pgf rec}
G_{i+1}(\mathbf{z})=\e[\prod_{j\neq i} z_j^{X_i^j} e^{-V_i \sum_{j\neq i} \lambda_j (1-z_j)} e^{-D_i \sum_{j=1}^N \lambda_j (1-z_j)} \big(p_i(V_i)+z_i[1-p_i(V_i)]\big)^{X_i^i} e^{-\Lambda_i(V_i)(1-z_i)}].
\end{equation}
Recall that for all $j$, the random variables $X_i^j$ are independent of $V_i$ and $D_i$. Therefore, Equation~\eqref{eq:pgf rec} becomes
\begin{multline}\label{eq:intermediate}
G_{i+1}(\mathbf{z})=\e[e^{-D_i \sum_{j=1}^N \lambda_j (1-z_j)}] \times \\ \times \e[\prod_{j\neq i} z_j^{X_i^j} e^{-V_i \sum_{j\neq i} \lambda_j (1-z_j)} \big(p_i(V_i)+z_i[1-p_i(V_i)]\big)^{X_i^i} e^{-\Lambda_i(V_i)(1-z_i)}].
\end{multline}
Consequently,
\begin{multline}\label{eq:pgf}
G_{i+1}(\mathbf{z})=\ltd[i]\bigl(\,\sum_{j=1}^N \lambda_j (1-z_j)\bigr)\times\\
\times \e[e^{-V_i \sum_{j\neq i} \lambda_j (1-z_j)} e^{-\Lambda_i(V_i)(1-z_i)}G_i(z_1,z_2,\ldots,z_{i-1},p_i(V_i)+[1-p_i(V_i)] z_i,z_{i+1},\ldots,z_N)],
\end{multline}
\noindent where $\ltd[i](s)=\e[e^{-sD_i}]$ denotes the Laplace-Stieltjes transform of the random variable $D_i$. Evidently, if $V_i$ follows a discrete distribution, the above expression simplifies significantly. Note that the mean queue length at a polling instant \eqref{eq:means} can also be obtained by differentiating Equation~\eqref{eq:pgf}.

\begin{remark}
Applying similar techniques, we can also derive the probability generating function of the number of customers at the end of a visit at queue $i+1$. If we denote by $Y_i^j$ the number of customers in queue $j=1,\ldots,N$ at the moment when the service at queue $i=1,\ldots,N$ is completed, then  the law of motion describing the evolution of the system is given by
\begin{equation}\label{eq:lom2}
  Y_{i+1}^j=\begin{cases}
                Y_i^j+A_j(D_i)+A_j(V_{i+1}),     &j\neq i+1,\\
                \mathrm{Binom}(Y_i^j+A_j(D_i), 1-p_j(V_j))+\mathrm{Poisson}\big(\Lambda_j(V_j)\big), &j=i+1.
  \end{cases}
\end{equation}
Also note that the expected value of $Y_i^j$ can be easily computed from \eqref{eq:means} by observing that for all $j\neq i$, $Y_i^j=X_i^j+A_j(V_i)$, while for the $i$-th queue we have that $Y_i^i=\mathrm{Binom}(X_i^i, 1-p_i(V_i))+\mathrm{Poisson}\bigl(\Lambda_i(V_i)\bigr)$.
\end{remark}

\section{Sojourn time}\label{s:sojourn}
Let the sojourn time of a customer at queue $i$ be denoted by $S_i$. We compute its expected value (and thus, by Little's law, also the mean queue length of queue $i$ at an arbitrary moment), and we derive the Laplace-Stieltjes transform of $S_i$. As stated before, for each queue we assume that if the service of a customer is not completed during a visit, then, for the next visit at that queue, a new service time will be resampled for the same customer from the service time distribution of that queue.

\subsection{Mean sojourn time}
Recall that the cycle time is given by $C=\sum_{i=1}^N (V_i+D_i)$. In order to derive the mean sojourn time of a customer arriving at queue $i$, we shall need some further notation. Denote by $V_i^{res}$ the residual visit time of the group of servers at queue $i$ and by $C_{/i}$ the cycle time except the time spent serving queue $i$, i.e.\ $C_{/i}=C-V_i$. Similarly, $C_{/i}^{res}$ represents the residual cycle time excluding the visit time of queue $i$. That is, $C_{/i}^{res}$ measures the length of time from a random moment after leaving queue $i$ until the next polling instant of queue $i$. Furthermore, let $\{C_m\}$ be a family of i.i.d.\ random variables distributed like $C$, and $N_i$ be a (shifted) geometric random variable with success probability $p_i={\e[p_i(V_i)]=}\p[B_i\leqslant V_i]$, i.e.\ $\p[N_i=n]=(1-p_i)^n p_i$, for all integer $n\geqslant 0$. One should observe here that $N_i +1$ is a stopping time as it is the first time when the service time of a customer in queue $i$ is less than or equal to the visit time at that queue; that is, $N_i+1=\inf\{k: B_{i,k}\leqslant V_{i,k}\}$, where $B_{i,k}$ and $V_{i,k}$ are i.i.d.\ copies of $B_{i}$ and $V_{i}$ respectively. Similarly, a second index is added to a random variable, every time that we explicitly need to indicate that an independent copy is considered. Then the sojourn time of a customer at queue $i$ is given by
\begin{equation}\label{eq:sojourn}
S_i=\begin{cases}
        B_{i,0}, &\mbox{(arrival during $V_i$ and $B_{i,0} \leqslant V_{i,0}^{res}$),}\\
        V_{i,0}^{res}+\sum_{m=1}^{N_i} C_m +C_{/i}+B_{i,N_i+1}, &\mbox{(arrival during $V_{i,0}$ and $B_{i,0} > V_{i,0}^{res}$),}\\
        C_{/i}^{res}+ \sum_{m=1}^{N_i} C_m +B_{i,N_i+1}, &\mbox{(arrival during $C_{/i}$)}.
    \end{cases}
\end{equation}
Note that the probability of an arrival occurring during the visit time of queue $i$ is $\e[V_i]/\e[C]$, i.e.\ the expected visit time of queue $i$ over the expected cycle time, and similarly for the other two events. Therefore, from \eqref{eq:sojourn} we obtain that the expected sojourn time of a customer of queue $i$ is given by
\begin{multline}\label{eq:mean sojourn1}
\e[S_i]=\frac{\e[V_i]}{\e[C]}\ \p[B_i \leqslant V_i^{res}]\ \e[B_i \mid B_{i} \leqslant V_{i}^{res}]+\\
+\frac{\e[V_i]}{\e[C]}\ \p[B_i>V_i^{res}]\ \e[V_{i,0}^{res}+\sum_{m=1}^{N_i} C_m +C_{/i}+B_{i,N_i+1}  \mid B_{i,0} \geqslant V_{i,0}^{res}]+\\
+\frac{\e[C_{/i}]}{\e[C]}\ \e[C_{/i}^{res}+ \sum_{m=1}^{N_i} C_m +B_{i,N_i+1}].
\end{multline}
In order to compute the second conditional expectation appearing at the right-hand side of the above equation, we think as follows. For $N_i$ cycles, the service of the customer is not completed during that visit because for every visit $B_i>V_i$, while at the $N_i+1$st visit the service is completed within that cycle. Therefore, define
$$
\overline{C}_m=C_{/i,m}+\min(B_{i,m},V_{i,m})
$$
and observe that
$$
\e[\sum_{m=1}^{N_i+1} \overline{C}_m]=\e[\sum_{m=1}^{N_i} {C_m}+C_{/i}+B_{i}].
$$
Thus
\begin{align}
\nonumber \e[V_{i,0}^{res}+\sum_{m=1}^{N_i} C_m +C_{/i}&+B_{i,N_i+1}  \mid B_{i,0} \geqslant V_{i,0}^{res}]
    =\e[V_{i}^{res}\mid B_{i} \geqslant V_{i}^{res}]+\e[\sum_{m=1}^{N_i+1} \overline{C}_m]\\
\nonumber    &=\e[V_{i}^{res}\mid B_{i} \geqslant V_{i}^{res}]+\e[N_i+1]\e[\overline{C}_m]\\
             &=\e[V_{i}^{res}\mid B_{i} \geqslant V_{i}^{res}]+\e[N_i+1]\big(\e[C_{/i}]+\e[\min(B_{i},V_{i})]\big),
\end{align}
where in the second equality we used Wald's equation.

For the third conditional expectation appearing at the right-hand side of \eqref{eq:mean sojourn1}, we have that
\begin{align}
\nonumber \e[C_{/i}^{res}+ \sum_{m=1}^{N_i} C_m +B_{i,N_i+1}]&=\frac{\e[C_{/i}^2]}{2\e[C_{/i}]}+\e[\sum_{m=1}^{N_i} (C_{/i,m}+V_{i,m})+B_{i,N_i+1}]\\
\nonumber    &=\frac{\e[C_{/i}^2]}{2\e[C_{/i}]}+\e[\sum_{m=1}^{N_i} C_{/i,m}]+\e[\sum_{m=1}^{N_i+1} \min(B_{i,m},V_{i,m})]\\
            &=\frac{\e[C_{/i}^2]}{2\e[C_{/i}]}+\e[{N_i}]\e[ C_{/i}]+\e[N_i+1]\e[ \min(B_{i},V_{i})].
\end{align}
Summarising the above, we have that
\begin{multline}\label{eq:final mean}
\e[S_i]=\frac{\e[V_i]}{\e[C]}\ \p[B_i \leqslant V_i^{res}]\ \e[B_i \mid B_{i} \leqslant V_{i}^{res}]+\\
+\frac{\e[V_i]}{\e[C]}\ \p[B_i>V_i^{res}]\Big(\e[V_{i}^{res}\mid B_{i} > V_{i}^{res}]+\e[N_i+1]\big(\e[C_{/i}]+\e[\min(B_{i},V_{i})]\big)\Big)+\\
+\frac{\e[C_{/i}]}{\e[C]} \Big(\frac{\e[C_{/i}^2]}{2\e[C_{/i}]}+\e[{N_i}]\e[ C_{/i}]+\e[N_i+1]\e[ \min(B_{i},V_{i})]\Big).
\end{multline}
In Section~\ref{s:Numerical results} we shall illustrate through an example the effect of the first two moments of the visit time and the service time on the mean sojourn time of an arbitrary customer.

\subsection{The Laplace-Stieltjes transform}
We now derive the Laplace-Stieltjes transform of the sojourn time of a customer of queue $i$. We first rewrite Equation~\eqref{eq:sojourn} in terms of the Laplace-Stieltjes transforms of all variables involved (cf.\ \eqref{eq:final mean}), and thus we get that
\begin{multline}\label{eq:LST sojourn prep}
\e[e^{-s S_i}]=\frac{\e[V_i]}{\e[C]}\ \p[B_i \leqslant V_i^{res}] \e[e^{-s B_i}\mid B_{i} \leqslant V_{i}^{res}]+\\
+\frac{\e[V_i]}{\e[C]}\ \p[B_i>V_i^{res}] \e[e^{-s V_i^{res}}\mid B_{i} > V_{i}^{res}] \e[e^{-s \sum_{m=1}^{N_i+1} \overline{C}_m}]+\\
+\frac{\e[C_{/i}]}{\e[C]}\ \e[e^{-s C_{/i}^{res}}] \e[e^{-s \sum_{m=1}^{N_i} C_{/i,m}}] \e[e^{-s \sum_{m=1}^{N_i+1} \min(B_{i,m},V_{i,m})}].
\end{multline}
We rewrite several of the terms appearing above as follows. The distribution function of $V_{i}^{res}$ is given by
$$
\p[V_{i}^{res} \leqslant x]=\frac{1}{\e[V_{i}]}\int_0^x \p[V_{i} >y ] \d y,
$$
yielding
\begin{equation}\label{eq:prob B ge Vk res}
\p[B_i > V_i^{res}]=\frac{1}{\e[V_i]}\int_0^\infty \p[B_i >x]\p[V_i>x] \d x.
\end{equation}
Similarly, we have that
$$
\p[C_{/i}^{res} \leqslant x]=\frac{1}{\e[C_{/i}]}\int_0^x \p[C_{/i} >y ] \d y,
$$
which implies that
\begin{equation}\label{eq:LT C/i}
  \e[e^{-s C_{/i}^{res}}]=\frac{1-\widetilde{C}_{/i} (s)}{s \e[C_{/i}]},
\end{equation}
where $\widetilde{C}_{/i}$ denotes the Laplace-Stieltjes transform of the random variable $C_{/i}$. Moreover,
\begin{align}
\nonumber \e[e^{-s \sum_{m=1}^{N_i+1} \overline{C}_m}]&=\sum_{n=0}^\infty \e[e^{-s \sum_{m=1}^{n+1} \overline{C}_m}] (1-p_i)^n p_i=\sum_{n=0}^\infty \e[e^{-s \overline{C}}]^{n+1} (1-p_i)^n p_i \\
&=\frac{p_i \e[e^{-s \overline{C}}]}{1-(1-p_i)\e[e^{-s \overline{C}}]},
\end{align}
where $\overline{C}=C_{/i}+\min(B_{i},V_{i})$. Likewise, we have that
\begin{equation}
\e[e^{-s \sum_{m=1}^{N_i} C_{/i,m}}]=\frac{p_i}{1-(1-p_i)\e[e^{-s C_{/i}}]}
\end{equation}
and
\begin{equation}\label{eq:teleytaia}
\e[e^{-s \sum_{m=1}^{N_i+1} \min(B_{i,m},V_{i,m})}]=\frac{p_i \e[e^{-s \min(B_{i},V_{i})}]}{1-(1-p_i)\e[e^{-s \min(B_{i},V_{i})}]}.
\end{equation}

Substituting \eqref{eq:prob B ge Vk res} -- \eqref{eq:teleytaia} into \eqref{eq:LST sojourn prep} we have that the Laplace-Stieltjes transform of the sojourn time of a customer of queue $i$ is given by
\begin{multline}\label{eq:LT S}
\e[e^{-s S_i}]=\frac{1}{\e[C]}\ \e[e^{-s B_i}\mid B_{i} \leqslant V_{i}^{res}] \int_0^\infty \p[B_i\leqslant x]\p[V_i>x] \d x +\\
+\frac{1}{\e[C]}\ \e[e^{-s V_i^{res}}\mid B_{i} > V_{i}^{res}] \frac{p_i \e[e^{-s \overline{C}}]}{1-(1-p_i)\e[e^{-s \overline{C}}]} \int_0^\infty \p[B_i >x]\p[V_i>x] \d x+\\
+\frac{1-\widetilde{C}_{/i} (s)}{s \e[C]}\ \frac{p_i}{1-(1-p_i)\widetilde{C}_{/i} (s)} \frac{p_i \e[e^{-s \min(B_{i},V_{i})}]}{1-(1-p_i)\e[e^{-s \min(B_{i},V_{i})}]}.
\end{multline}
Clearly, from the expression above, one can retrieve Equation~\eqref{eq:final mean} for the mean sojourn time of a customer of queue $i$.

The transforms appearing in \eqref{eq:LT S} may be cumbersome to compute when the service times or the visit times are generally distributed. However, when both $B_i$ and $V_i$ follow a phase-type distribution, all transforms can be computed explicitly since the class of phase-type distributions is closed under finite minima. Phase-type distributions are widely used in computations. The class of phase-type distributions is dense (in the sense of weak convergence) in the class of all distributions on $(0,\infty)$  (cf.\ \cite[Propositions 1 and 2]{asmussen00a}). As an example, we will derive the Laplace-Stieltjes transform of the sojourn time of a customer of queue $i$, as well as its mean, in case both the visit time and the service time at queue $i$ are exponentially distributed.

\subsection{A special case}\label{ss:{A special case}}
Let the service time and the visit time at queue $i$ be exponentially distributed with rates $\mu_i$ and $\gamma_i$ respectively. Then all terms appearing in \eqref{eq:final mean} can be easily computed in terms of $\mu_i$ and $\gamma_i$. For example,
$$
\e[B_i \mid B_{i} \leqslant V_{i}^{res}]=\frac{1}{\gamma_i +\mu_i }
$$
and $\p[B_i>V_i^{res}]=\gamma_i /(\gamma_i +\mu_i )$. Thus, \eqref{eq:final mean} becomes
\begin{multline*}
\e[S_i]=\frac{1}{\gamma_i \e[C]}\ \frac{\mu_i}{\gamma_i +\mu_i }\ \frac{1}{\gamma_i +\mu_i }
+\frac{1}{\gamma_i \e[C]}\ \frac{\gamma_i}{\gamma_i +\mu_i }\Big(\frac{1}{\gamma_i +\mu_i }+\big(\frac{\gamma_i}{\mu_i}+1\big)\big(\e[C_{/i}]+\frac{1}{\gamma_i +\mu_i }\big)\Big)+\\
+\frac{\e[C_{/i}]}{\e[C]} \Big(\frac{\e[C_{/i}^2]}{2\e[C_{/i}]}+\frac{\gamma_i}{\mu_i}\e[ C_{/i}]+(\frac{\gamma_i}{\mu_i}+1)\frac{1}{\gamma_i +\mu_i }\Big)
\end{multline*}
or
$$
\e[S_i]=\frac{(\gamma_i \e[C_{/i}]+1)^2}{\gamma_i \mu_i \e[C]}+\frac{\e[C_{/i}^2]}{2\e[C]}.
$$

Similarly, \eqref{eq:LT S} reduces to
\begin{multline*}
\e[e^{-s S_i}]=\frac{1}{\e[C]}\ \frac{\gamma_i +\mu_i}{\gamma_i +\mu_i+s} \frac{\mu_i}{\gamma_i (\gamma_i+\mu_i)}
+\frac{1}{\e[C]}\ \frac{\gamma_i +\mu_i}{\gamma_i +\mu_i+s} \frac{\mu_i \e[e^{-s \overline{C}}]}{\gamma_i+\mu_i-\gamma_i\e[e^{-s \overline{C}}]} \frac{1}{\gamma_i+\mu_i}+\\
+\frac{1-\widetilde{C}_{/i} (s)}{s \e[C]}\ \frac{\mu_i}{\gamma_i+\mu_i-\gamma_i\widetilde{C}_{/i} (s)} \frac{\mu_i \frac{\gamma_i +\mu_i}{\gamma_i +\mu_i+s}}{\gamma_i+\mu_i-\gamma_i\frac{\gamma_i +\mu_i}{\gamma_i +\mu_i+s}}.
\end{multline*}
Since $\e[e^{-s \overline{C}}]=\widetilde{C}_{/i} (s)\e[e^{-s \min(B_{i},V_{i})}]$ we have that the previous expression reduces to
\begin{multline*}
\e[e^{-s S_i}]=\frac{1}{\e[C]}\ \frac{1}{\gamma_i +\mu_i+s} \frac{\mu_i}{\gamma_i }
+\frac{1}{\e[C]}\ \frac{1}{\gamma_i +\mu_i+s} \frac{\mu_i \widetilde{C}_{/i} (s)}{\gamma_i+\mu_i+ s -\gamma_i\widetilde{C}_{/i} (s)} +\\
+\frac{1-\widetilde{C}_{/i} (s)}{s \e[C]}\ \frac{\mu_i}{\gamma_i+\mu_i-\gamma_i\widetilde{C}_{/i} (s)} \frac{\mu_i} {\mu_i+s}.
\end{multline*}
Similar expressions can be easily derived in case both the visit times and the service times follow some phase-type distribution, such as Gamma, hyperexponential, or Coxian distributions.

\section{Numerical results}\label{s:Numerical results}
This section is devoted to some numerical results. In particular, we want to examine numerically the effect of the first two moments of the visit and service times on the sojourn time of an arbitrary customer. In all examples, we make the following assumptions. We consider a polling system with two queues. The arrival rate at the first queue is $\lambda_1=0.8$ and at the second queue it is $\lambda_2=0.5$. The service time and the visit time at the first queue are exponentially distributed with rates $\mu_1=1$ and $\gamma_1=1$ respectively. Moreover, the total mean switch-over time is taken to be $\e[D]=0.5$, while its second moment is assumed to be zero. In all figures that follow, we plot the mean sojourn time of an arbitrary customer, which is estimated by $(\lambda_1 \e[S_1]+\lambda_2 \e[S_2])/(\lambda_1+\lambda_2)$.

In Figures~\ref{fig:ESvsEB} and \ref{fig:ESvsCB} we investigate the effect of the first two moments of the service time at the second queue on the mean sojourn time of an arbitrary customer. For these plots, the visit time at the second queue is considered to be exponentially distributed with rate $\gamma_2=3/2$. For various values of the squared coefficient of variation of the service time at the second queue, which is denoted by $c_{B_2}^2$, we plot in Figure~\ref{fig:ESvsEB} the mean sojourn time of an arbitrary customer versus the mean service time $\e[B_2]$. The squared coefficient of variation of the service time is chosen to be comparable to the squared coefficient of variation of the (exponentially distributed) visit time, which is equal to 1. In Figure~\ref{fig:ESvsCB}, we plot the mean sojourn time of an arbitrary customer versus $c_{B_2}^2$ for three values of $\e[B_2]$, which again are chosen to be comparable to $\e[V_2]$.

For each case of $c_{B_2}^2$, we fit a mixed Erlang or hyperexponential distribution to $\e[B_2]$ and $c_{B_2}^2$, depending on whether the squared coefficient of variation is less or greater than one; see, e.g., Tijms~\cite{tijms-FCSM}. So, if $1/n \leqslant c_{B_2}^2 \leqslant 1/(n-1)$ for some $n = 2, 3, \ldots$, then the mean and squared coefficient of variation of the mixed-Erlang distribution
$$
G(x)= p\,
\biggl(1-\mathrm{e}^{-\zeta x}\sum_{j=0}^{n-2}\frac{(\zeta x)^j}{j!}\biggr)
+ (1-p)
\biggl(1-\mathrm{e}^{-\zeta x}\sum_{j=0}^{n-1}\frac{(\zeta x)^j}{j!}\biggr),
\qquad x \geqslant 0 ,
$$
matches with $\e[B_2]$ and $c_{B_2}^2$, provided the parameters $p$ and $\zeta$ are chosen as
$$
p = \frac{1}{1+c_{B_2}^2}\Bigl(n c_{B_2}^2-\sqrt{n(1+c_{B_2}^2)-n^2 c_{B_2}^2}\ \Bigr),\qquad\zeta=\frac{n-p}{\e[B_2]}.
$$
On the other hand, if $c_{B_2}^2>1$, then the mean and squared coefficient of variation of the hyperexponential distribution
$$
G(x) = p (1- \mathrm{e}^{-\zeta_1 x})+q(1-\mathrm{e}^{-\zeta_2 x}),\qquad x \geqslant 0,
$$
match with $\e[B_2]$ and $c_{B_2}^2$, provided the parameters $\zeta_1$, $\zeta_2$, $p$, and $q$ are chosen as
\begin{align*}
&p =\frac{1}{2}\biggl(1+\sqrt{\frac{c_{B_2}^2-1}{c_{B_2}^2+1}}\,\biggr),  &&q=1-p,\\
&\zeta_1 =\frac{2 p}{\e[B_2]} &\mbox{and}\quad&\zeta_2=\frac{2 q}{\e[B_2]}.
\end{align*}

As is evident from the plot in Figure~\ref{fig:ESvsEB}, the expected sojourn time of an arbitrary customer increases as the mean service time at the second queue increases. Moreover, the rate that it increases with is almost linear as $c_{B_2}^2$ grows and the effect of the second moment is less pronounced than the effect of $\e[B_2]$.

In Figure~\ref{fig:ESvsCB}, one observes that the mean sojourn time of an arbitrary customer decreases as the squared coefficient of variation of the service time increases, contrary to what is the case for the M/G/1 queue. This result is due to the fact that the service time of a customer that did not complete his service during one visit time is resampled for the following visit time. Therefore, the larger the variability in the service times, the bigger is the probability that during the next visit time this particular customer will complete his service. Recall that the mean visit time at the second queue is equal to 2/3 and observe that in case $\e[B_2]$ is less than $\e[V_2]$, the effect of the second moment of the service time on the mean sojourn time of an arbitrary customer is almost negligible.

\noindent\begin{minipage}[t]{0.49\textwidth}
\includegraphics[width=\textwidth]{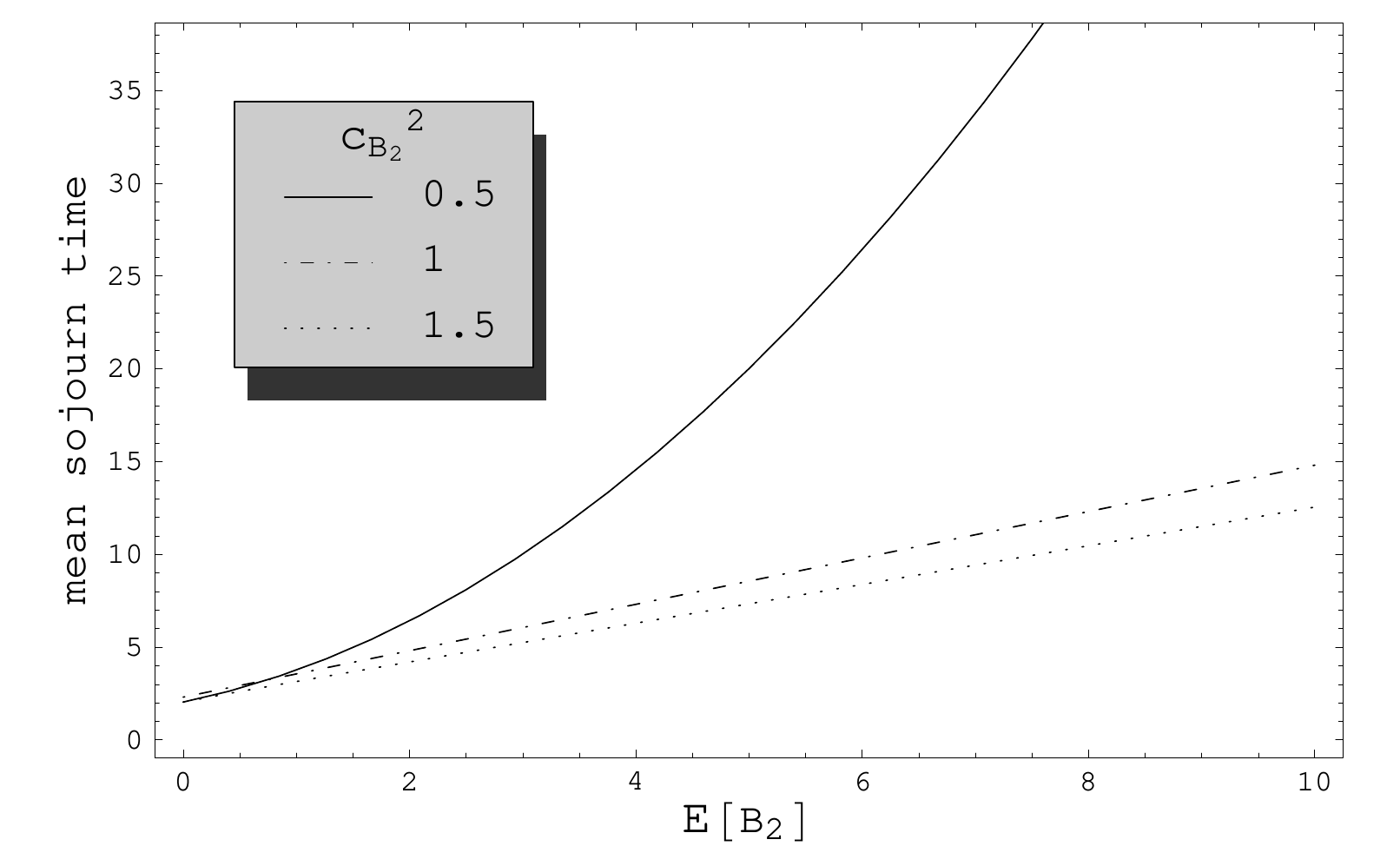}
\figcaption{Mean sojourn time of an arbitrary customer against the mean service time $\e[B_2]$.}
\label{fig:ESvsEB}
\end{minipage}
\hfill
\begin{minipage}[t]{0.49\textwidth}
\includegraphics[width=\textwidth]{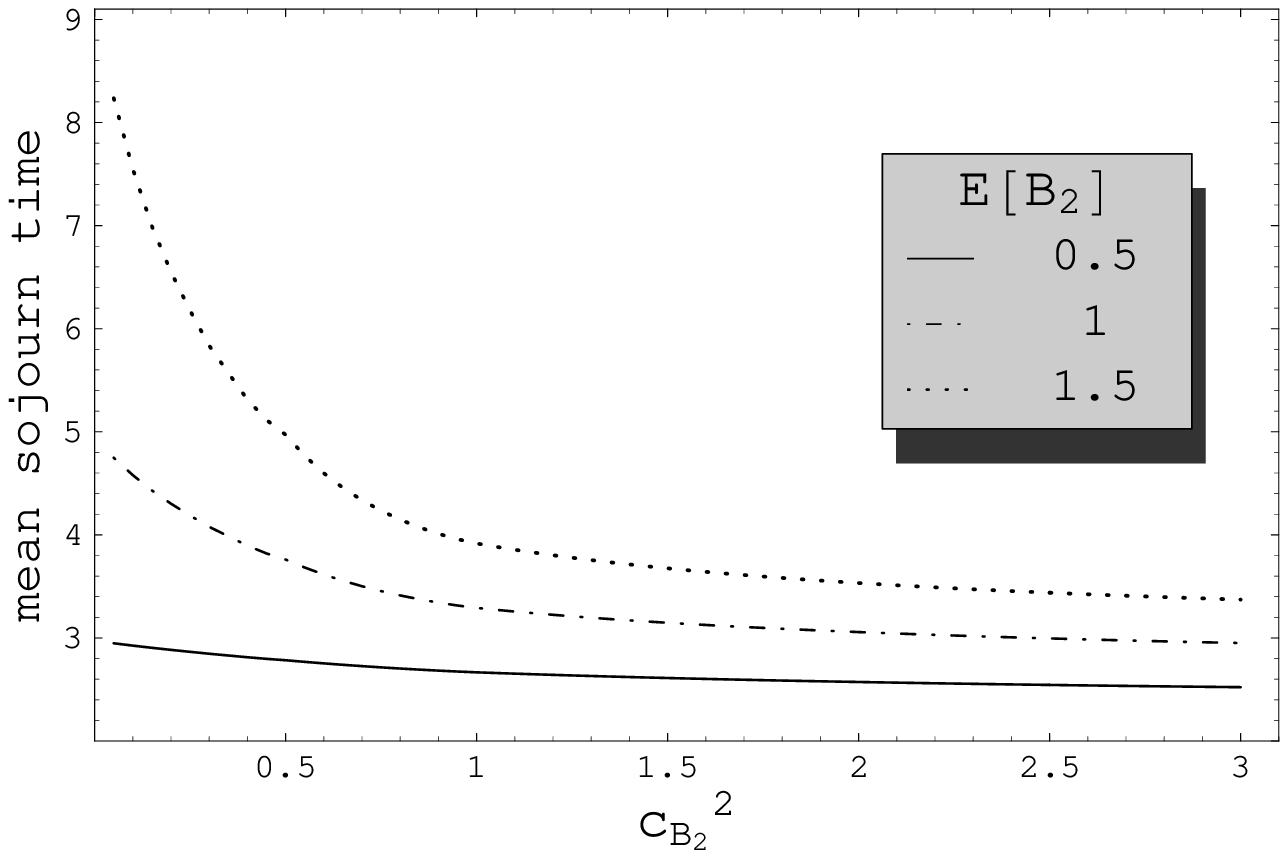}
\figcaption{Mean sojourn time of an arbitrary customer against the squared coefficient of variation of the service time $B_2$.}
\label{fig:ESvsCB}
\end{minipage}
\vspace*{\baselineskip}

In Figures~\ref{fig:ESvsEV} and \ref{fig:ESvsCV} we now investigate the effects of the first two moments of the visit time at the second queue on the mean sojourn time of an arbitrary customer. For these plots, we now take the service time at the second queue to be exponentially distributed with rate $\mu_2=3/2$. For various values of the squared coefficient of variation of the visit time at the second queue, which is denoted by $c_{V_2}^2$, we plot in Figure~\ref{fig:ESvsEV} the mean sojourn time of an arbitrary customer versus the mean visit time $\e[V_2]$. As before, the squared coefficient of variation of the visit time is chosen to be comparable with the squared coefficient of variation of the (exponentially distributed) service time, which is equal to 1. In Figure~\ref{fig:ESvsCV}, we plot the mean sojourn time of an arbitrary customer versus $c_{V_2}^2$ for three values of $\e[V_2]$, which again are chosen to be comparable with $\e[B_2]$.

The plot in Figure~\ref{fig:ESvsEV} is interesting. Evidently, when $\e[V_2]$ is significantly smaller than $\e[B_2]$, only a very small number of customers will be served during a visit. As the mean visit time increases, more customers are served during a visit and the mean sojourn time of an arbitrary customer is reduced. However, as the mean visit time continues to increase, this trend is reversed after the mean sojourn time of an arbitrary customer reaches a global minimum. In other words, there is an optimal value for the mean visit time to some queue that minimises the mean sojourn time of an arbitrary customer; beyond that value, the mean sojourn time of an arbitrary customer increases at an almost linear rate. This indicates that the polling system under consideration can be optimised in expectation by controlling the visit time to each queue. In the following section, we will develop a policy that minimises the mean sojourn time of an arbitrary customer in the system.

\begin{minipage}[t]{0.47\textwidth}
\includegraphics[width=\textwidth]{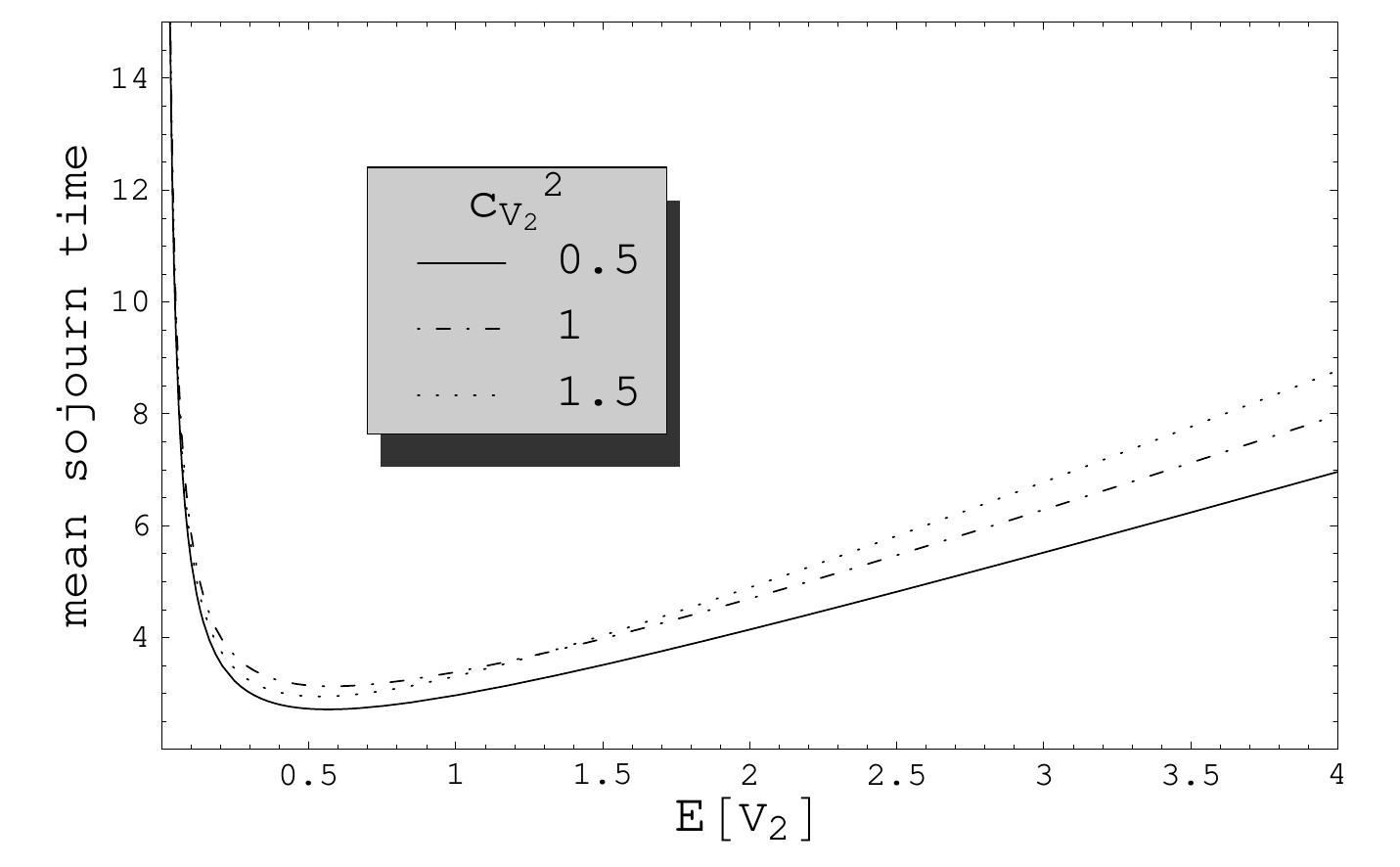}
\figcaption{Mean sojourn time of an arbitrary customer against the mean visit time $\e[V_2]$.}
\label{fig:ESvsEV}
\end{minipage}
\hfill
\begin{minipage}[t]{0.47\textwidth}
\includegraphics[width=\textwidth]{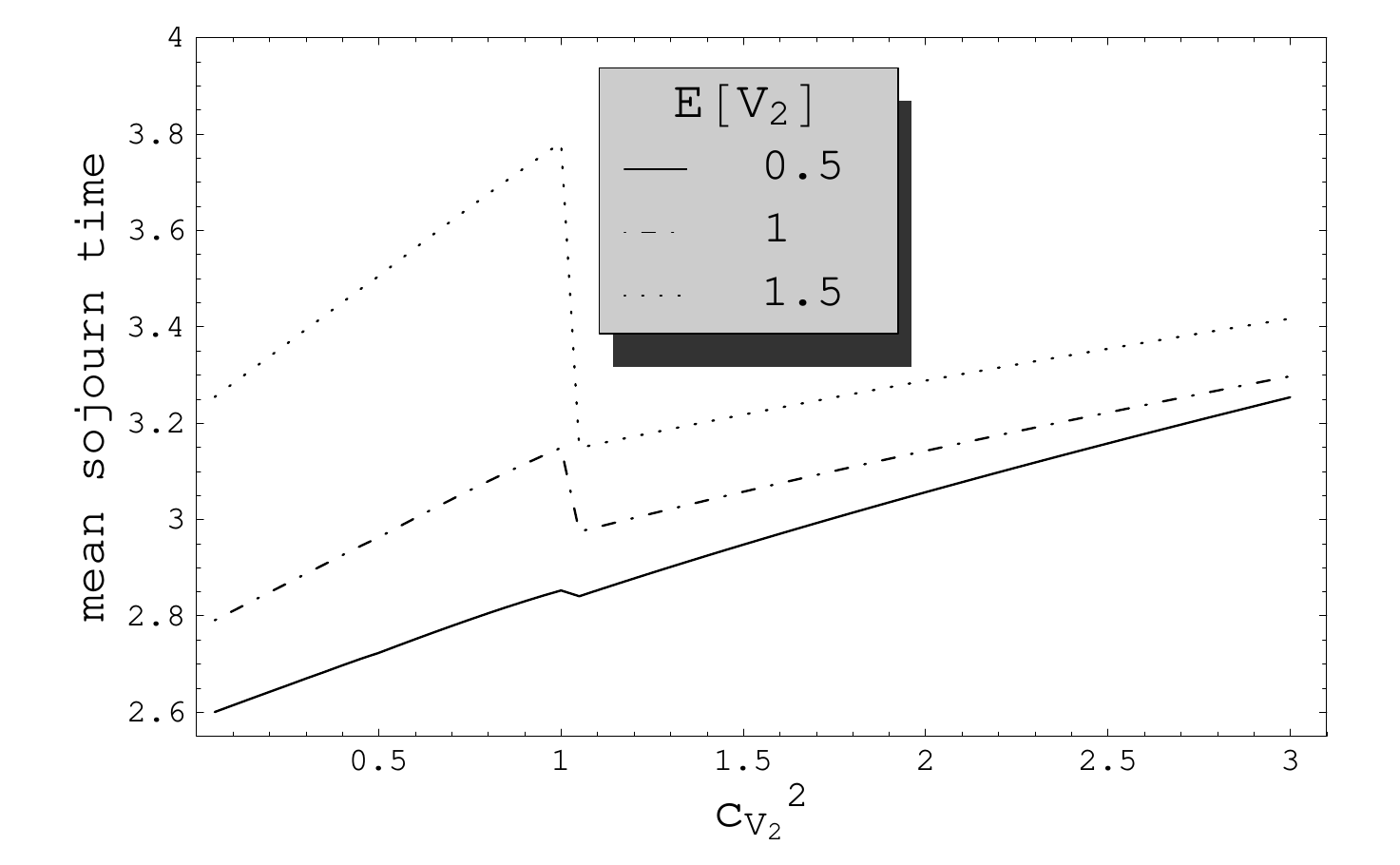}
\figcaption{Mean sojourn time of an arbitrary customer against the squared coefficient of variation of the visit time $V_2$.}
\label{fig:ESvsCV}
\end{minipage}
\vspace*{\baselineskip}

The plot in Figure~\ref{fig:ESvsCV} is also interesting. As is explained above, we fit either a mixture of Erlang distributions or a hyperexponential distribution to each pair of the first two moments of the visit time, depending on the value of the squared coefficient of variation. For every value of $c_{V_2}^2$, we obtain a different visit time distribution. Note that the jump in Figure~\ref{fig:ESvsCV} occurs when the distribution we fit to the first two moments of the visit time shifts from a mixture of Erlang distributions to a hyperexponential distribution. This indicates that the shape of the visit time distribution is important; for example, hyperexponential distributions are always unimodal, which is not the case for mixed Erlang distributions. Consequently, the first two moments cannot capture sufficiently the effect of the visit time distribution on the mean sojourn time of an arbitrary customer; one needs to know the exact distribution.

\section{Dynamic control of servers' visits}\label{s:optimisation}
A basic question that arises when planning efficient polling systems concerns the order of visits performed by the servers.  As it is suggested by Figure~\ref{fig:ESvsEV}, the polling system we are considering can be optimised in some way so that the mean sojourn time of an arbitrary customer is minimised. Rather than identifying the value of the minimum mean sojourn time for the cyclic processing (visiting) order considered so far, we first investigate whether there exists a fixed static order that the servers visit the various queues so that the mean (weighted) sojourn time of an arbitrary customer is minimised. As is evident from Equation~\eqref{eq:final mean}, the mean sojourn time of a customer of queue $i$ \emph{does not depend} on the order the queues have been visited, and thus neither does the weighted sum thereof. Since the mean sojourn time of an arbitrary customer remains unaffected when altering the processing order, it does not constitute a practical performance measure of our system.

An appealing approach that leads to a simple and tractable rule is to develop a semi-dynamic control scheme. The idea is to dispatch the group of servers to perform Hamiltonian tours, each tour being possibly different from the previous one, depending on the state of the system at the beginning of the tour, so as to optimise some performance measure. An adequate performance measure is the \emph{throughput} of the system, namely the number of customers served per cycle, as the throughput can be measured per cycle, while the sojourn time of a customer spans over a random number of cycles. The goal is to maximise the throughput of the system for each cycle.

Specifically, suppose that at the beginning of a cycle $n=(n_1, n_2,\ldots, n_N)$ is the state of the system, where $n_i$ is the number of jobs waiting in queue $i$. We shall compute the expected value of the number of customers served per cycle under a specific processing order, and consequently identify the optimal processing order per cycle. The following theorem summarises the result.

\begin{theorem}
For the Hamiltonian-tour approach, the optimal visiting order is \emph{independent} of the number of customers present at the various queues at the start of the cycle and is given by the index-type rule
$$
\frac{\lambda_i p_i}{\e[V_i]+\e[D_i]}$$
in the sense that the throughput is maximised if and only if the visiting order is arranged according to an \emph{increasing} sequence of this rule.
\end{theorem}

Before proceeding with the proof we point out that having more information regarding the system, such as the number of customers present at all queues, has no effect on the optimal strategy, and thus it does not improve the performance of the system. This stands in contrast to many other polling systems where typically more information regarding a system leads to decisions that increase the efficiency of the system; see for example \cite{browne88}, \cite{browne89}, and \cite{yechiali93}. This conclusion stems from the fact that no matter what cycle order is used, the number of customers served in a cycle among those initially present will have the same distribution and will have no effect on the number of others served (or the number receiving partial service) during the cycle. In other words, what happens to those initially present at the beginning of a tour is unaffected by what ordering is used.

\begin{proof}
The proof follows from an interchange argument. Consider the processing order $\pi_0=(1,2,\ldots,N)$. Denote by $\theta_i$ the throughput of queue $i$ under this processing order, i.e.\ the number of customers of queue $i$ that are served during a cycle, and denote by $\theta$ the total throughput of the system, i.e.\ the sum of all $\theta_i$. Given $n$, i.e.\ the state of the system at the beginning of a cycle, we shall compute the expected value of $\theta$.

The number of customers served after completing a visit at queue $i$ is equal to the portion of customers that where present at polling instant $i$ and successfully completed their service plus the number of customers that arrived during the visit time of queue $i$ and completed their service within that visit. In other words,
\begin{equation*}
\theta_i=\mathrm{Binom}\Bigl(n_i+A_i\bigl(\,\sum_{k=1}^{i-1} (V_k+D_k)\bigr), p_i(V_i)\Bigr)+ \Bigl(A_i(V_i)-\mathrm{Poisson}(\Lambda_i(V_i)\Bigr).
\end{equation*}
As a result,
$$
\e[\theta_i]=\Bigl(n_i+\lambda_i \sum_{k=1}^{i-1}\bigl( \e[V_k]+\e[D_k]\bigr)\Bigr) p_i +\lambda_i \e[V_i]-\e[\Lambda_i(V_i)],
$$
which yields
\begin{equation}\label{eq:e-theta}
\e[\theta]=c+\sum_{i=1}^N \lambda_i p_i \sum_{k=1}^{i-1} \bigl(\e[V_k]+\e[D_k]\bigr),
\end{equation}
where
$$
c=\sum_{i=1}^N (n_i p_i + \lambda_i \e[V_i]-\e[\Lambda_i(V_i)]).
$$
Observe that the constant $c$ that appears in \eqref{eq:e-theta} does not depend on $\pi_0$, while the second term at the right-hand side of \eqref{eq:e-theta} does.

Consider now the processing order $\pi_1=(1,2,\ldots,j-1,j+1,j,j+2,\ldots,N)$, where the visit order of queues $j$ and $j+1$ is interchanged and denote by $\theta^{\,\prime}_i$ and $\theta^{\,\prime}$ the throughput of queue $i$ and of the whole system under $\pi_1$, respectively. We promptly have that $\e[\theta_i]=\e[\theta^{\,\prime}_i]$ for all $i\neq j, j+1$ and that
\begin{align*}
&  \e[\theta^{\,\prime}_j]=\Bigl(n_j+\lambda_j \bigl(\,\sum_{k=1}^{j-1} (\e[V_k]+\e[D_k])+\e[V_{j+1}]+ \e[D_{j+1}]\bigr)\Bigr) p_j +\lambda_j \e[V_j]-\e[\Lambda_j(V_j)],\\
&  \e[\theta^{\,\prime}_{j+1}]=\Bigl(n_{j+1}+\lambda_{j+1} \sum_{k=1}^{j-1}\bigl( \e[V_k]+\e[D_k]\bigr)\Bigr) p_{j+1} +\lambda_{j+1} \e[V_{j+1}]-\e[\Lambda_{j+1}(V_{j+1})].
\end{align*}
Thus,
\begin{multline*}
\e[\theta^{\,\prime}]=c+\sum_{i\neq j, j+1} \lambda_i p_i \sum_{k=1}^{i-1} \bigl(\e[V_k]+\e[D_k]\bigr)+\lambda_j p_j \bigl(\,\sum_{k=1}^{j-1} (\e[V_k]+\e[D_k])+\e[V_{j+1}]+ \e[D_{j+1}]\bigr)+\\
+\lambda_{j+1} p_{j+1} \sum_{k=1}^{j-1}\bigl(\e[V_k]+\e[D_k]\bigr).
\end{multline*}
Therefore, we have that $\e[\theta]\leqslant\e[\theta^{\,\prime}]$ if and only if
\begin{multline*}
\lambda_j p_j \sum_{k=1}^{j-1} \bigl(\e[V_k]+\e[D_k]\bigr)+\lambda_{j+1} p_{j+1} \sum_{k=1}^{j} \bigl(\e[V_k]+\e[D_k]\bigr)\leqslant\\
\lambda_j p_j \sum_{k=1}^{j-1} \bigl(\e[V_k]+\e[D_k]\bigr)+\lambda_j p_j \bigl(\e[V_{j+1}]+ \e[D_{j+1}]\bigr)+\lambda_{j+1} p_{j+1} \sum_{k=1}^{j-1}\bigl(\e[V_k]+\e[D_k]\bigr),
\end{multline*}
or
$$
\lambda_{j+1} p_{j+1} \bigl(\e[V_j]+\e[D_j]\bigr)\leqslant\lambda_j p_j \bigl(\e[V_{j+1}]+ \e[D_{j+1}]\bigr).
$$
In other words, we get that the optimal processing order is by visiting the queues according to an \emph{increasing} order of $\lambda_i p_i/\bigl(\e[V_i]+\e[D_i]\bigr)$.
\end{proof}

Roughly stated, this rule  arranges the visit order according to the ratio between new arrivals per unit time that will successfully complete their service, i.e.\ $\lambda_i p_i$, and the mean duration of a visit there, i.e.\ $\e[V_i]+\e[D_i]$. It is intuitively clear that if the mean visit and switch time for a queue is relatively long, then one should visit this queue early on. In this way, the number of customers at the other queues during this cycle will also be relatively high, and as a result the throughput will be increased since all customers are served simultaneously by an infinite number of servers.

This is an extremely simple rule, which can be directly implemented. Moreover, suppose that, for one reason or another, the objective is to minimise the throughput of the system for each cycle. Then, the index rule that determines the order of visits to the queues is simply \emph{reversed}; the servers complete a Hamiltonian tour arranged in a decreasing order of $\lambda_i p_i/\bigl(\e[V_i]+\e[D_i]\bigr)$. Observe that under this strategy the servers also visit the queues that are empty at the beginning of the cycle.

One expects that the throughput of the system in the long-run is improved when these queues are not visited within a cycle; namely, it may be more efficient to avoid queues that are empty at the beginning of the cycle in order to allow them to build up.

According to the way the system is designed, even if the servers do not visit a queue that at the beginning of the cycle was empty, the switch time associated with this queue (i.e.\ the time to switch from this queue to the following one) is still incurred. Therefore, as the number of queues that will not be visited in a cycle grows, the servers spend an increasing amount of time being essentially idle (as they switch between queues).

A possibly more efficient system design is as follows. Rather than envision the group of servers moving from one queue to another, we can think of a central point to which the servers always return after each completion of a visit to a queue. The return time to the central point after visiting queue $i$ is denoted by $R_i$. The servers depart from that central point and move to the following queue that will be served. The total time from the moment the servers leave the central point until they enter queue $i$ is denoted by $E_i$. According to this design, the total time to go from queue $i$ to queue $j$ is given by $R_i+E_j$ for any $i\neq j$. The question that arises is whether there exists a semi-dynamic control of this system. As before, it emerges that a Hamiltonian-tour approach leads to a \emph{static} processing order according to an index rule.

\begin{theorem}
For the polling system with a central point, the Hamiltonian-tour approach leads to a fixed optimal visiting order, which is \emph{independent} of the number of customers present at the various queues at the start of the cycle. The throughput of the system for each cycle is maximised if and only if the visiting order is arranged according to an \emph{increasing} sequence of the index-type rule
\begin{equation}
\label{IndexRule2}
\frac{\lambda_i p_i}{\e[E_i]+\e[V_i]+\e[R_i]}.
\end{equation}
\end{theorem}
\begin{proof}
As before, let $n=(n_1, n_2,\ldots, n_N)$ be the state of the system at the start of the tour and denote by $L$ the number of \emph{non-empty} queues, $0<L\leqslant N$, at the beginning of the cycle. The throughput of queue $i$ during a Hamiltonian cycle that visits only the non-empty queues according to the order $\pi_0=(1,2,\ldots,L)$ is given by
\[
\theta_i=\mathrm{Binom}\Bigl(n_i+A_i\bigl(\,\sum_{k=1}^{i-1} (E_k+V_k+R_k)+E_i\bigr), p_i(V_i)\Bigr)+ \Bigl(A_i(V_i)-\mathrm{Poisson}(\Lambda_i(V_i))\Bigr).
\]
Consequently,
$$
\e[\theta_i]=\biggl(n_i+\lambda_i \Bigl(\sum_{k=1}^{i-1}\bigl(\e[E_k]+\e[V_k]+\e[R_k]\bigr)+\e[E_i]\Bigr)\biggr) p_i +\lambda_i \e[V_i]-\e[\Lambda_i(V_i)],
$$
which yields
\begin{equation}\label{eq:e-thetanew}
\e[\theta]=c^\prime+\sum_{i=1}^N \lambda_i p_i \sum_{k=1}^{i-1} \bigl(\e[E_k]+\e[V_k]+\e[R_k]\bigr),
\end{equation}
where
$$
c^\prime=\sum_{i=1}^N \bigl((n_i + \lambda_i \e[E_i])p_i + \lambda_i \e[V_i]-\e[\Lambda_i(V_i)]\bigr).
$$
Applying an interchange argument we have that the optimal processing order is constructed by an increasing sequence of the index rule given by \eqref{IndexRule2}.
\end{proof}

As before, the optimal tour does not depend on the number of customers present at the beginning of the cycle. This is a direct consequence of the fixed visit times and the underlying M/G/$\infty$ process at each queue.

Index rules appear regularly when optimising polling systems. Browne and Yechiali~\cite{browne88,browne89} were the first to obtain dynamic control policies for single-server systems under the exhaustive, gated or mixed service regimes. The mechanics of the system are as described here: at the beginning of each cycle the server decides on a new Hamiltonian tour and visits the channels accordingly. The authors showed that if the objective is to optimise the cycle duration under these policies, then an index-type rule applies, which is similar to the one described here. The main difference is that the index rule that is optimal for these policies depends on the state of the system at the beginning of a cycle, \emph{contrary to the results obtained for the fixed-visit-time policy} studied in this paper. The result derived by Browne and Yechiali~\cite{browne88,browne89} is a surprising result as the index rule does \emph{not} include the service times at the various channels. It is also surprising that the same index rule holds for both the gated and the exhaustive disciplines although the duration of a cycle starting from the same state is \emph{different} for the two regimes. For a further discussion on other types of index-rule policies see Yechiali~\cite{yechiali93}, van der Wal and Yechiali~\cite{wal03}, and references therein.

\section*{Acknowledgements}
The authors would like to thank Prof.\ Onno Boxma for several insightful remarks. The authors also acknowledge the hospitality and support of EURANDOM while carrying out this research. The research of the first author was also supported by the Aristotle University of Thessaloniki (full scholarship from the legacy of L.\ Athanasoula). Part of the second author's work was carried out while visiting EURANDOM in his capacity as Beta Chair.

%%\bibliographystyle{abbrv}
%%\bibliography{maria}

\end{document}